\documentclass[12pt]{article}

\usepackage{fullpage}
\usepackage{amssymb}
\usepackage{amsmath}
\usepackage{amsthm}
\usepackage{url}
\usepackage{hyperref}

\newtheorem{theorem}{Theorem}

\newtheorem{lemma}{Lemma}

\begin{document}

\title{Permutations $r_j$ such that $\sum_{i=1}^n \prod_{j=1}^k r_j(i)$ is maximized or minimized}

\date{August 12, 2015\\Latest update: February 29, 2020}
\author{Chai Wah Wu\\ IBM T. J. Watson Research Center\\ P. O. Box 218, Yorktown Heights, New York 10598, USA\\e-mail: chaiwahwu@ieee.org}
\maketitle

\begin{abstract}
We consider the problem of finding the set of permutations $r_j$ of $\{1,\cdots , n\}$ such that $\sum_{i=1}^n \prod_{j=1}^k r_j(i)$ is maximized or minimized. While the set of permutations maximizing this value are easily determined, finding the set of permutations minimizing this value appears to be an open problem. We show values of $k$ and $n$ for which an explicit solution exists and comment on computational issues in determining the general problem. We also look at the dual problem of finding the permutations such that $\prod_{i=1}^n \sum_{j=1}^k r_j(i)$ is maximized or minimized. As part of this study we also look at a variant of a rearrangement inequality.
\end{abstract}

\section{Introduction}\label{sec:intro}
Let $0 \leq a_1 \leq  a_2  \dots \leq a_n$ be a sequence of nonnegative numbers. 
Consider $k$ permutations of the integers $\{1,\cdots, n\}$ denoted as $\{r_1,\cdots, r_k\}$
and define the value 
$v(n,k) = \sum_{i=1}^n \prod_{j=1}^k a_{r_j(i)}$.  The maximal and minimal value of $v$ among all $k$-sets of permutations are denoted as $v_{\max}(n,k)$ and $v_{\min}(n,k)$ respectively.
For $k=2$, $v$ corresponds to the inner product between two permutations of the sequence $\{a_i\}$. The goal of this paper is to study the values of $v_{\max}$ and $v_{\min}$ and the sets of permutations that achieve them.

It turns out the maximal value $v_{\max}$ is easy to determine whereas $v_{\min}$ is not. In particular, it is easy to show that $v_{\max}(n,k) = \sum_{i=1}^n a_i^k$ and 
is achieved when 
all the $k$ permutations are the same.  
This is a consequence of the following result in \cite{ruderman}:

\begin{lemma}\label{lem:ruderman}
Consider a set of nonnegative numbers $\{a_{ij}\}$, $i=1,\cdots, k$, $j=1,\cdots, n$.  Let
$a'_{i1},a'_{i2},\cdots,a'_{in}$ be the numbers $a_{i1},a_{i2},\cdots,a_{in}$ reordered such that
$a'_{i1}\geq a'_{i2}\geq\cdots\geq a'_{in}$.  Then
\[ \sum_{j=1}^n\prod_{i=1}^k a_{ij} \leq   \sum_{j=1}^n\prod_{i=1}^k a'_{ij} \]
\[ \prod_{j=1}^n\sum_{i=1}^k a_{ij} \geq   \prod_{j=1}^n\sum_{i=1}^k a'_{ij} \]
\end{lemma}

Recall the AM-GM inequality:
\begin{lemma}\label{lem:am-gm}
For $n$ nonnegative real numbers $x_i\geq 0$, $\sum_{i=1}^n x_i \geq n \sqrt[n]{\prod_{i=1}^n x_i}$ and   $\prod_{i=1}^n x_i \leq \left(\frac{\sum_{i=1}^n x_i}{n}\right)^n$ with equality if and only if all the $x_i$ are the same.
\end{lemma}
The use of AM-GM inequality allows us to get a lower bound on $v_{\min}$.
\begin{lemma} \label{lem:vmin}
$v_{\min}(n,k) \geq  n \prod_i a_i^{k/n}$.
\end{lemma}
\begin{proof}
The product $\prod_{ij} a_{r_i(j)}$ is equal to $\prod_i a_i^k$. Thus by Lemma \ref{lem:am-gm}, $v(n,k) \geq n\sqrt[n]{ \prod_i a_i^k} = n \prod_i a_i^{k/n}$.
\end{proof}

The cases of $n=1$ or $k=1$ are simple as clearly $v(1,k) = a_1^k$ and $v(n,1) = \sum_{i=1}^n a_i$.

\section{The case $n=2$}
There are only two permutations on the integers $\{1,2\}$.
In this case $v_{\max}(2,k) = a_1^k+a_2^k$.  If $k=2m$ is even, $v_{\min}(2,k) = 2a_1^ma_2^m$ is achieved with $k/2$ of the permutations of one kind and the other half the other kind.
If $k=2m+1$ is odd, $v_{\min}(2,k) = (a_1+a_2)a_1^ma_2^m$ is achieved with $m$ of the permutations of one kind and $m+1$ of them the other kind.

\section{The case $k=2$}

\begin{lemma}[Rearrangement inequality]\label{lem:rearrange}
\[x_ny_1+\cdots +x_1y_n \leq x_{\sigma(1)}y_1+\cdots + x_{\sigma(n)}y_n \leq x_1y_1+\cdots x_ny_n\]

for real numbers $x_i$, $y_i$ such that
$x_1\leq \cdots \leq x_n$ and $y_1\leq \cdots \leq y_n$ and all permutations $\sigma$.
\end{lemma}
A proof of this can be found in \cite{hardy}.

Using this, we can prove the following:

\begin{theorem} $v_{\max}(n,2) = \sum_{i=1}^n a_i^2$ and  $v_{\min}(n,2) = \sum_{i=1}^n a_ia_{n-i+1}$.
\end{theorem}
\begin{proof}
This is clear from Lemma \ref{lem:rearrange} by choosing both permutations to be  ($1$,$2$,$\cdots$, $n$) for $v_{\max}(n,2)$
and choosing the two permutations to be ($1$,$2$,$\cdots$, $n$) and ($n$,$n-1$,$\cdots$ , $2$, $1$) $v_{\min}(n,2)$.  
\end{proof}

Finding $v_{\min}(n,k)$ appears to be nontrivial for $k > 2$ and $n > 2$.
Since $v(n,k)$ is invariant under simultaneous reordering of the permutations, we can use this to define equivalence classes among $k$-sets of permutations.

\section{The case where $k$ is a multiple of $n$}
Rudeman's inequality (Lemma \ref{lem:ruderman}) is an extension of Lemma \ref{lem:rearrange} and its dual Lemma \ref{lem:dual_rearrange} (Section \ref{sec:dual}) to multiple permutations. However, whereas Lemmas \ref{lem:rearrange} and \ref{lem:dual_rearrange} provide both an upper and lower bound, Lemma \ref{lem:ruderman} only provides an upper bound for the sum of products and a lower bound for the product of sums. In this section we provides conditions for which a lower bound and an upper bound can be obtained for the sum of products and product of sums respectively. As a result of proving these conditions we also provide a characterization of the permutation that achieves $v_{\min}(n,k)$, when $n$ divides $k$.

\begin{theorem}
If $n$ divides $k$, then $v_{min}(n,k) = \left(\prod_{i=1}^n a_i^{k/n}\right) n$ and is achieved by using each cyclic permutation $k/n$ times. 
\end{theorem}

\begin{proof}
By Lemma \ref{lem:vmin} $v(n,k) \geq  \left(\prod_{i=1}^n a_i^{k/n}\right) n$.
Consider the $n$ cyclic permutations $\sigma_1 = (1,2,...,n)$, $\sigma_2 = (2,...,n,1)$, ..., $\sigma_n = (n,1,...,n-1)$.
It is clear that using $k/n$ copies of each permutation $\sigma_i$ to form $k$ permutations results in $v =\left(\prod_{i=1}^n a_i^{k/n}\right) n$.
\end{proof}

\section{Computing $v_{\min}(n,k)$}
Since we can always pick the member in the equivalent class such that $r_1 = \{1,\cdots, n\}$, we can fix $r_1$.  By Lemma \ref{lem:rearrange}, after $r_1,\cdots r_{k-1}$ are chosen, the permutation $r_k$ that minimize
$v(n,k)$ among all choices of $r_k$ is the permutation of $\{1,\cdots, n\}$ that is in reverse order from the order of the sequence of numbers $\left\{\prod_{j=1}^{k-1} a_{r_j(1)}, \prod_{j=1}^{k-1} a_{r_j(2)}, \cdots , \prod_{j=1}^{k-1} a_{r_j(n)}\right\}$.  This means we only need to test $v(n,k)$ by choosing only $r_2,\cdots, r_{k-1}$, as $r_1$ is fixed and $r_k$ is determined by the choices for the other permutations.
The number of combinations we need to check (for $k > 2$) is equal to the number of combinations of $n!$ objects chosen $k-2$ times with replacement, which is equal to 
\[\left(\begin{array}{c}n!+k-3\\k-2\end{array}\right)\]

The following Python function computes $v_{\min}(n,k)$ (for $k \geq 2$):

\begin{verbatim}
from itertools import permutations, combinations_with_replacement
def vmin(n,k,a): # compute vmin(n,k) for k > 1
    a = sorted(a)
    ntuple, count = tuple(range(n)), sum(d**k for d in a)
    for s in combinations_with_replacement(permutations(ntuple,n),k-2):
        t = list(a)
        for d in s:
            for i in range(n):
                t[i] *= a[d[i]]
        t.sort()
        v = 0
        for i in range(n):
            v += a[n-i-1]*t[i]
            if v >= count:
                break
        if v < count:
            count = v
    return count        
\end{verbatim}

Further reduction in the number of permutations searched can be obtained for large $k$. Since the number of permutations is $n!$, given $k$ permutations, by the extended pigeonhole principle, there will be a permutation that is repeated $t = \lfloor\frac{k-1}{n!}\rfloor + 1$ times. Thus we can fix
$r_1, \cdots, r_t$ to be the same permutation $\{1,\cdots n\}$ and choose only $r_{t+1}\cdots, r_{k-1}$. Thus the number of combinations we need to check is
\[\left(\begin{array}{c}n!+k-t-2\\k-t-1\end{array}\right)\]

\section{Special case where $a_i = i$}
Consider the special case where the sequence of numbers $\{a_i\}$ are the integers $\{1,\cdots , n\}$, i.e. $a_i = i$.
In this case $v_{\max}(n,k) = \sum_{i=1}^n i^k$.
When $k$ is much smaller than $n$, computation of this value can be simplified via Faulhaber's formula as:
\[ v_{\max}(n,k) = \frac{n^{k+1}}{k+1} + \frac{1}{2}n^k + \sum_{i =2}^{k} \frac{B_i}{i!}\frac{k!}{(k-i+1)!}n^{k-i+1}\]
where $B_j$ are the Bernouilli numbers.

$v_{\max}(n,n) = \sum_{i=1}^n i^n = H(n,-n)$, the $n$-th generalized harmonic number of order $-n$ (OEIS sequence A031971, \url{https://oeis.org/A031971}) 

Next, we have $v(1,k) = 1$,$v(n,1) = \sum_{i=1}^n i = \frac{n(n+1)}{2}$. Furthermore, 
if $n$ divides $k$, then $v_{min}(n,k) = \left(\prod_{i=1}^n i^{k/n}\right) n = n!^{k/n}n$. This implies that $v_{\min}(n,n) = n!n = (n+1)!-n!$ (OEIS sequence A001563, \url{https://oeis.org/A001563}).

Furthermore, $v_{\max}(2,k) = 1+2^k$.  If $k=2m$ is even, $v_{\min}(2,k) = 2^{m+1}$
If $k=2m+1$ is odd, $v_{\min}(2,k) = 3\cdot 2^m$.

$v_{\max}(n,2) = \sum_{i=1}^n i^2 = \frac{n(n+1)(2n+1)}{6}$.
Consider the two permutations ($1$,$2$,$\cdots$, $n$) and ($n$,$n-1$,$\cdots$ , $2$, $1$).  The value of $v(n,2)$ is 
equal to $\sum_{i=1}^n i(n-i+1) = (n+1)\sum_{i=1}^n i - \sum_{i=1}^n i^2 = \frac{n(n+1)^2}{2} - \frac{n(n+1)(2n+1)}{6} = 
n(n+1)\left(\frac{n+1}{2}-\frac{2n+1}{6}\right) = \frac{n(n+1)(n+2)}{6}$ and is in fact equal to $v_{\min}(n,2)$ by Lemma \ref{lem:rearrange}.

\subsection{$v_{\min}(n,k)$ for $k=3,4, \cdots$}
The values of $v_{\min}(n,k)$ for different $k$'s are listed in OEIS sequences \cite{oeis} A070735 ($k=3$, \url{https://oeis.org/A070735}), A070736 ($k=4$, \url{https://oeis.org/A070736}), A260356 ($k=5$, \url{https://oeis.org/A260356}), A260357 ($k=6$, \url{https://oeis.org/A260357}), A260358 ($k=7$, \url{https://oeis.org/A260358}), and sequence A260359 (for the case $k=8$, \url{https://oeis.org/A260359}).

Partial list of values of $v_{\min}(n,k)$ (with some data taken from OEIS) are listed in Table \ref{tbl:vmin}. 
The antidiagonals of Table \ref{tbl:vmin} can be found in \url{https://oeis.org/A260355} (OEIS sequence A260355).

\begin{table}[htbp]
\begin{center}
\tiny
\begin{tabular}{r||r|r|r|r|r|r|r|r|r|r|r|r|r|}
& $k=1$ & 2 & 3&4&5&6&7&8&9&10&11&12&13\\
\hline\hline
$n=1$&1    & 1&     1  &   1 &    1  &   1 &    1  &   1   &  1 &    1  &   1   &  1  &    1  \\ \hline
2& 3  &   4 &    6 &    8   & 12  &  16 &   24  &  32 &   48  &  64 &   96 &  128 &   192   \\ \hline
3 &6 &   10  &  18 &   33  &  60 &  108 &  198  & 360 &  648 & 1188 & 2145 & 3888  & 7083  \\ \hline
4&10&    20   & 44 &   96  & 214 &  472 & 1043 & 2304 & 5136 &11328 &24993 &55296& 122624   \\ \hline
5 & 15   & 35  &  89 &  231 &  600 & 1564 & 4074 &10618&27665&72000&&&\\\hline
6 & 21  &  56  & 162 &  484 & 1443 & 4320 &&&&&&3110400&\\\hline
7 &  28 &   84 &  271  & 915 & 3089 &&35280&&&&&&\\\hline
8 & 36 &  120 &  428 & 1608 &&&&322560&&&&& \\\hline
9 & 45 &  165 &  642 & 2664 &&&&&3265920&&&& \\\hline
10 & 55  & 220 &  930 & 4208 &&&&&&3628800&&& \\\hline
11 & 66 &  286 & 1304 &&&&&&&&439084800&& \\\hline
12 & 78  & 364 & 1781 &&&&&&&&&5748019200& \\\hline
13 & 91  & 455 & 2377 &&&&&&&&&&80951270400 \\\hline
14 & 105  & 560 & 3111 &&&&&&&&&& \\\hline
15 & 120  & 680 & 4002 &&&&&&&&&& \\\hline
16 & 136  & 816 & 5073 &&&&&&&&&& \\\hline
17 & 153 & 969 & 6344 &&&&&&&&&& \\\hline
18 & 171  & 1140 & 7842 &&&&&&&&&& \\\hline
19 & 190  & 1330 & 9587 &&&&&&&&&& \\\hline
\end{tabular}
\end{center}

\begin{center}
\tiny
\begin{tabular}{r||r|r|r|r|r|r|r|r|r|r|r|r|}
& $k=14$& 15&16 &17 &18 &19&20&21&22&23&24&25\\
\hline\hline
$n=1$&1&1 &1    & 1&     1  &   1 &    1  &   1 &    1  &   1   &  1 &    1  \\ \hline
2& 256&384 & 512  &   768 &   1024 &    1536   & 2048  &  3072 &   4096  &  6144 &   8192  &  12288  \\ \hline
3 &12844&23328 &42498 &  77064  & 139968   & 254988& 462384    & 839808   & 1526769  & 2774304   & 5038848  & 9160614   \\ \hline
4&271040&599832 & 1327104 &       &  &     &31850496  &   &  &  &764411904  &  \\ \hline
\end{tabular}
\end{center}

\begin{center}
\tiny
\begin{tabular}{r||r|r|r|r|r|r|r|r|r|r|r|r|r|r|r|r|}
& $k=26$&27&28&29&30\\
\hline\hline
$n=1$  & 1   &  1  &    1  &    1  &    1\\ \hline
2&    16384 &  24576 &   32768  &  49152  &  65536 \\ \hline
3 &16645824 &30233088  &54963684  &   99874944& 181398528     \\ \hline
4&  &  &18345885696 && \\ \hline
\end{tabular}
\end{center}

\caption{Partial list of $v_{\min}(n,k)$.}\label{tbl:vmin}
\end{table}

Table \ref{tbl:vmin2} is a partial table of the number of nonequivalent $k$-sets of permutations achieving $v_{\min}(n,k)$ where equivalence is described in Section \ref{sec:intro}.  We will denote these numbers as $N_{\min}(n,k)$.
For $n \leq 2$ or $k\leq 2$, $N_{\min}(n,k) = 1$.  On the other hand, $N_{\min}(n,k)$ can be larger than $1$ if $n >2$ or $k > 2$.
For example, the $2$ sets of nonequivalent permutations that achieves $v_{\min}(3,6) = 108$ are
$(123, 123, 231, 231, 312, 312)$ and
$(123, 132, 213, 231, 312, 321)$.
The $3$ sets of nonequivalent permutations that achieves $v_{\min}(5,3) = 89$ are
$(12345, 34251, 52314)$, 
$(12345, 35214, 52341)$ and
$(12345, 35241, 52314)$.

Note that $N_{\max}(n,k)$ corresponding to $v_{max}$ satisfies $N_{\max}(n,k) = 1$ for all $n$ and $k$.

\begin{table}[htbp]
\begin{center}
\footnotesize
\begin{tabular}{r||r|r|r|r|r|r|r|r|r|r|r|r|r|r|r|}
& $k=1$ & 2 & 3&4&5&6&7&8&9&10&11&12&13&14&15\\
\hline\hline
$n=1$&1    & 1&     1  &   1 &    1  &   1 &    1  &   1   &  1 &    1  &   1   &  1  &    1  &    1  &    1  \\ \hline
2&1    & 1&     1  &   1 &    1  &   1 &    1  &   1   &  1 &    1  &   1   &  1  &    1  &    1  &    1  \\ \hline
3 &1  &   1  &   1  &   1  &   1  &   2 &    1  &   2  &   2  &   2  &   1 &    3 &    1   &   1  &    3 \\ \hline
4& 1  &   1  &   2  &   4  &  11  &  10  &  10  &  81  & 791  & 533 & 24 & 1461 & 3634 & 192 & 2404\\ \hline
5 & 1  &   1   &  3 &   12 &   16 &  188 &  211& 2685&&&&&&&\\\hline
6 & 1  &   1  &  10 &  110 &   16             &&&&&&&&&&\\\hline
7 &  1  &   1  &  6 &   &                   &&&&&&&&&&\\\hline
8 & 1  &   1  &  16 &   &                   &&&&&&&&&&\\\hline
9 &  1  &   1  &  4 &   &                   &&&&&&&&&&\\\hline
10 & 1  &   1  &  12 &   &                    &&&&&&&&&&\\\hline
11 &  1  &   1  &   &   &                   &&&&&&&&&&\\\hline
12 &  1  &   1  &   &   &                   &&&&&&&&&&\\\hline
13 &  1  &   1  &   &   &                   &&&&&&&&&&\\\hline
14 & 1  &   1  &   &   &                    &&&&&&&&&&\\\hline
15 & 1  &   1  &   &   &                    &&&&&&&&&&\\\hline
\end{tabular}
\end{center}
\caption{Partial list of $N_{\min}(n,k)$, the number of nonequivalent $k$-sets of permutations that achieve $v_{\min}(n,k)$.}.\label{tbl:vmin2}
\end{table}

Table \ref{tbl:Nmin} lists $N_{\min}(3,k)$ for various values of $k$.

\begin{table}[htbp]
\begin{center}
\begin{tabular}{r|r}
$k$ & $N_{\min}(3,k)$ \\
\hline\hline
1 & 1\\ \hline
2 & 1\\ \hline
3 & 1\\ \hline
4 & 1\\ \hline
5 & 1\\ \hline
6 & 2\\ \hline
7 & 1\\ \hline
8 & 2\\ \hline
9 & 2\\ \hline
10 & 2\\ \hline
11 & 1\\ \hline
12 & 3\\ \hline
13 & 1\\ \hline
14 & 1\\ \hline
15 & 3\\ \hline
16 & 2\\ \hline
17 & 1\\ \hline
18 & 4\\ \hline
19 & 3\\ \hline
20 & 2\\ \hline
21 & 4\\ \hline
22 & 1\\ \hline
23 & 2\\ \hline
24 & 5\\ \hline
25 & 1\\ \hline
26 & 3\\ \hline
27 & 5\\ \hline
28 & 2\\ \hline
29 & 3\\ \hline
30 & 6\\ \hline
31 & 2\\ \hline
32 & 4\\ \hline
33 & 6\\ \hline
34 & 3\\ \hline
35 & 1\\ \hline
\end{tabular}
\end{center}
\caption{$N_{\min}(3,k)$ for various values of $k$.} \label{tbl:Nmin}
\end{table}

The following Python function computes $v_{\min}(n,k)$ (for $k \geq 2$):

\begin{verbatim}
from itertools import permutations, combinations_with_replacement
def vmin(n,k): # compute vmin(n,k) for k > 1
    ntuple, count = tuple(range(1,n+1)), n**(k+1)
    for s in combinations_with_replacement(permutations(ntuple,n),k-2):
        t = list(ntuple)
        for d in s:
            for i in range(n):
                t[i] *= d[i]
        t.sort()
        v = 0
        for i in range(n):
            v += (n-i)*t[i]
            if v >= count:
                break
        if v < count:
            count = v
    return count
\end{verbatim}

\section{Permutation sets achieving $v_{\min}(n,k)$}
The next sections list, for each $n$ and $k$, a $k$-set of permutations (there may be many) that achieves $v_{\min}(n,k)$.  Concatenating the permutations and considering them as a number in base $n+1$, the listed $k$-set of permutations is the $k$-set with the smallest such number that achieves $v_{\min}(n,k)$.    We use the letters $a$, $b$, etc. to represent the numbers $10$, $11, \cdots $.
We omit the $n\leq2$ or $k\leq 2$ cases as they were discussed above.

\subsection{$k=3$}
\begin{itemize}
\item $n=3$: $(123, 231, 312)$
\item $n=4$: $(1234, 2341, 4213)$
\item $n=5$: $(12345, 34251, 52314)$
\item $n=6$: $(123456, 435261, 642315)$
\item $n=7$: $(1234567, 5463271, 7523416)$
\item $n=8$: $(12345678, 64572381, 86425317)$
\item $n=9$: $(123456789, 854673291, 976324518)$
\item $n=10$: $(123456789a, 96485372a1, a783452619)$
\item $n=11$: $(123456789ab, a65847932b1, b984632571a)$
\item $n=12$: $(123456789abc, b768593a24c1, ca856372941b)$ 
\end{itemize}

\subsection{$k=4$}
\begin{itemize}
\item $n=3$: $(123, 132, 312, 321)$
\item $n=4$: $(1234, 2143, 3412, 4321)$
\item $n=5$: $(12345, 23145, 42531, 54312)$
\item $n=6$: $(123456, 235146, 632541, 653412)$
\item $n=7$: $(1234567, 3264571, 6724513, 7542136)$
\item $n=8$: $(12345678, 32457168, 87423651, 87452613)$
\end{itemize}

\subsection{$k=5$}
\begin{itemize}
\item $n=3$: $(123, 123, 231, 312, 321)$
\item $n=4$: $(1234, 1234, 3214, 4231, 4321)$
\item $n=5$: $(12345, 21453, 34512, 45231, 53124)$
\item $n=6$: $(123456, 213564, 453612, 563241, 643125)$
\item $n=7$: $(1234567, 2317564, 5371624, 6574312, 7534162)$
\end{itemize}

\subsection{$k=6$}
\begin{itemize}
\item $n=3$: $(123, 123, 231, 231, 312, 312)$
\item $n=4$: $(1234, 1234, 2134, 4312, 4321, 4321)$
\item $n=5$: $(12345, 12345, 31254, 45213, 54321, 54321)$
\item $n=6$: $(123456, 123465, 421536, 564312, 635142, 654321)$
\end{itemize}

\subsection{$k=7$}
\begin{itemize}
\item $n=3$: $(123, 123, 132, 231, 312, 312, 321)$
\item $n=4$: $(1234, 1234, 1234, 4231, 4231, 4312, 4312)$
\item $n=5$: $(12345, 12345, 21534, 45132, 45231, 52314, 54321)$
\end{itemize}

\subsection{$k=8$}
\begin{itemize}
\item $n=3$: $(123, 123, 123, 231, 231, 312, 312, 321)$
\item $n=4$: $(1234, 1234, 1243, 3124, 3421, 4213, 4321, 4321)$
\item $n=5$: $(12345, 12345, 12345, 42513, 45123, 53142, 53421, 53421)$
\end{itemize}

\subsection{$k=9$}
\begin{itemize}
\item $n=3$: $(123, 123, 123, 231, 231, 231, 312, 312, 312)$
\item $n=4$: $(1234, 1234, 1234, 2134, 3241, 3412, 4213, 4321, 4321)$
\item $n=5$: $(12345, 12345, 12345, 42135, 42513,43251,43521,45213,54321)$
\end{itemize}

\subsection{$k=10$}
\begin{itemize}
\item $n=3$: $(123, 123, 123, 132, 231, 231, 312, 312, 312, 321)$
\item $n=4$: $(1234, 1234, 1234, 1234, 3124, 4213, 4321, 4321, 4321, 4321)$
\end{itemize}

\subsection{$k=11$}
\begin{itemize}
\item $n=3$: $(123, 123, 123, 123, 132, 312, 312, 321, 321, 321, 321)$
\item $n=4$: $(1234, 1234, 1234, 1234, 2314, 3412, 4213, 4231, 4231, 4231, 4231)$
\end{itemize}

\subsection{$k=12$}
\begin{itemize}
\item $n=3$: $(123, 123, 123, 123, 231, 231, 231, 231, 312, 312, 312, 312)$
\item $n=4$: $(1234, 1234, 1234, 1243, 2143, 3124, 3412, 3421, 4213, 4321, 4321, 4321)$

\end{itemize}

\subsection{$k=13$}
\begin{itemize}
\item $n=3$: $(123, 123, 123, 123, 132, 132, 312, 312, 312, 321, 321, 321, 321)$
\item $n=4$: $(1234, 1234, 1234, 1234, 1234, 2143, 4213, 4213, 4321, 4321, 4321, 4321, 4321)$
\end{itemize}

\subsection{$k=14$}
\begin{itemize}
\item $n=3$: $(123, 123, 123, 123, 123, 123, 213, 312, 321, 321, 321, 321, 321, 321)$
\item $n=4$: $(1234, 1234, 1234, 1234, 1234, 1324, 4132, 4132, 4312, 4312, 4321, 4321, 4321, 4321)$
\end{itemize}

\subsection{$k=15$}
\begin{itemize}
\item $n=3$: $(123, 123, 123, 123, 123, 231, 231, 231, 231, 231, 312, 312, 312, 312, 312)$
\item $n=4$: ($1234$, $1234$, $1234$, $1234$, $1234$, $1234$, $3214$, $3421$, $4213$, $4213$, $4231$, $4231$, $4231$, $4321$, $4321$)
\end{itemize}

\subsection{$k=16$}
\begin{itemize}
\item $n=3$: ($123$, $123$, $123$, $123$, $123$, $132$, $132$, $231$, $312$, $312$, $312$, $312$, $321$, $321$, $321$, $321$)
\item $n=4$: ($1234$, $1234$, $1234$, $1234$, $1243$, $1243$, $3124$, $3124$, $3421$, $3421$, $4213$, $4213$, $4321$, $4321$, $4321$, $4321$)
\end{itemize}

\subsection{$k=17$}
\begin{itemize}
\item $n=3$: ($123$, $123$, $123$, $123$, $123$, $123$, $213$, $231$, $312$, $312$, $321$, $321$, $321$, $321$, $321$, $321$, $321$)
\end{itemize}

\section{Permutations $r_j$ such that $\prod_{i=1}^n \sum_{j=1}^k a_{r_j(i)}$ is maximized or minimized} \label{sec:dual}
Consider the dual problem of determining the value 
$w(n,k) = \prod_{i=1}^n \sum_{j=1}^k a_{r_j(i)}$ and let the maximum and minimal value of $w$ among all $k$-sets of permutations be denoted as $w_{\max}(n,k)$ and $w_{\min}(n,k)$ respectively. The analysis of $w_{\max}$ and $w_{\min}$ is analogous to the analysis of $v_{\min}$ and $v_{\max}$, respectively.

By Lemma \ref{lem:ruderman},
$w_{\min}(n,k) = \prod_{i=1}^n k a_i = k^n \prod_i a_i $ obtained by choosing all the permutations $r_i$ to be the same.
In particular, $w_{\min}(2,k) = k^2 \prod_i a_i$ and $w_{\min}(n,2) = 2^n \prod_i a_i$.

Similarly, the AM-GM inequality (Lemma \ref{lem:am-gm}) gives the the following upper bound for $w_{\max}$.
\begin{lemma} \label{lem:ubwmax}
$w_{\max}(n,k) \leq \left(\frac{k\sum_i a_i}{n}\right)^n$.
\end{lemma}
\begin{proof}
The sum $\sum_{ij} r_i(j)$ is equal to $k\sum_i a_i$. Thus by Lemma \ref{lem:am-gm}, $w(n,k) \leq \left(\frac{k\sum_i a_i}{n}\right)^n$.
\end{proof}

On the other hand, similar to $v_{\min}(k,n)$, we have the following results for $w_{\max}(k,n)$ that are easy to prove:
\begin{itemize}
\item $w_{\max}(1,k) = ka_1$
\item $w_{\max}(n,1) = \prod_i a_i$
\item If $k = 2m$ is even, then $w_{\max}(2,k) = (a_1+a_2)^2m^2$.
\item If $k = 2m+1$ is odd, then $w_{\max}(2,k) = (ma_1+(m+1)a_2)(ma_2+(m+1)a_1)$.
\end{itemize}

There exists a product of sums version of the rearrangment inequality (Lemma \ref{lem:rearrange}):

\begin{lemma}[Dual rearrangement inequality]\label{lem:dual_rearrange}
\[(x_1+y_1)\times \cdots  (x_n+y_n) \leq (x_{\sigma(1)}+y_1)\times \cdots  (x_{\sigma(n)}+y_n ) \leq (x_n+y_1)\times \cdots (x_1+y_n)\]

for non-negative real numbers $x_i$, $y_i$ such that
$x_1\leq \cdots \leq x_n$ and $y_1\leq \cdots \leq y_n$ and all permutations $\sigma$.
\end{lemma}
This result first appeared in \cite{oppenheim}. See also \cite{geretschlager,ho}.
Note that unlike Lemma \ref{lem:rearrange}, the numbers $x_i$ and $y_i$ are required to be non-negative as Lemma \ref{lem:dual_rearrange} is not true in general for negative numbers.

\begin{theorem}
$w_{\max}(n,2) = \prod_i (a_i+a_{n-i+1})$.
\end{theorem}

\begin{theorem}
If $n$ divides $k$, then 
$w_{\max}(n,k) = \left(\frac{k\sum_i a_i}{n}\right)^n$ and is achieved by choosing each cyclic permutation $k/n$ times.
\end{theorem} 

In particular, $w_{\max}(n,n) = \left(\sum_i a_i\right)^n$.
The corresponding Python function for computing $w_{\max}(n,k)$ (for $k\geq 2$):

\begin{verbatim}
from itertools import permutations, combinations_with_replacement
def wmax(n,k,a): # compute wmax(n,k) for k > 1
    a = sorted(a)
    ntuple, count = tuple(range(n)), 0
    for s in combinations_with_replacement(permutations(ntuple,n),k-2):
        t = list(a)
        for d in s:
            for i in range(n):
                t[i] += a[d[i]]
        t.sort()
        w = 1
        for i in range(n):
            w *= a[n-i-1]+t[i]
        if w > count:
            count = w
    return count     
\end{verbatim}

The following table show various values of $w_{\max}(n,k)$.

\subsection{Special case of $a_i = i$}

The value $w_{\min}(n,k) = \prod_{i=1}^n ik = n!k^n$ is  obtained by choosing all the permutations $r_i$ to be the same and follows from Lemma \ref{lem:ruderman}.
$w_{\min}(n,n) = n!n^n$ (OEIS sequence A061711 \url{https://oeis.org/A061711}).

On the other hand, we have the following analogous results for $w_{\max}$:

\begin{theorem}
$w_{\max}$ satisfies the following equations:
\begin{itemize}
\item $w_{\max}(n,k) \leq \left(\frac{k(n+1)}{2}\right)^n$.
\item $w_{\max}(1,k) = k$
\item $w_{\max}(n,1) = n!$
\item If $k = 2m$ is even, then $w_{\max}(2,k) = 9m^2$.
\item If $k = 2m+1$ is odd, then $w_{\max}(2,k) = (3m+1)(3m+2)$.
\end{itemize}
\end{theorem}

\begin{theorem}
$w_{\max}(n,2) = (n+1)^n$.
\end{theorem}
In other words, $w_{\max}(n,2)$ is the number of labeled rooted trees with $n+1$ nodes (OEIS sequence A000169 \url{https://oeis.org/A000169}).

\begin{theorem}
If $n$ divides $k$, then 
$w_{\max}(n,k) = \left(\frac{k(n+1)}{2}\right)^n$.
\end{theorem} 
In particular, $w_{\max}(n,n) = \left(\frac{n(n+1)}{2}\right)^n$ (OEIS sequence A061718 \url{https://oeis.org/A061718}).
Interestingly, unlike $v_{\min}$, the value of $w_{\max}$ can be found for many values of $k$ and $n$.

\begin{theorem}
 $w_{\max}(3,k) = 8k^3$ for $k > 1$.
If $k$ is even, then $w_{\max}(3,k)$ is obtained by $k/2$ copies of the permutation $(1,2,3)$ followed by $k/2$ copies of the permutation $(3,2,1)$.
If $k = 2m+1\geq 3$ is odd, then $w_{\max}(3,k)$ is obtained by $m$ copies of the permutation $(1,2,3)$, the permutations $(2,3,1)$ and $(3,1,2)$ followed by $m-1$ copies of the permutation $(3,2,1)$.
\end{theorem}

\begin{proof}
By Lemma \ref{lem:ubwmax}, $w_{\max}(3,k) \leq 8k^3$. For both cases ($k$ is even and $k > 2$ is odd), the listed sequence of permutations shows that $w_{\max}(3,k) \geq 8k^3$.
\end{proof}

The proof of the following is straight forward.
\begin{theorem}
If $k$ is even, then $w_{\max}(n,k) = \left(\frac{k(n+1)}{2}\right)^n$ and is achieved by $k/2$ copies of the permutation 
$(1,\cdots, n)$ and $k/2$ copies of the permutation $(n,\cdots, 1)$.
\end{theorem}

\begin{proof}
The theorem follows from the fact that $\sum_{j}r_j(i) = k(n+1)/2$ for all $i$.
\end{proof}

The proof of the following result can be found in \cite{wu:rearrangement:2020}.
\begin{theorem}
If  $n$ is odd and $k\geq n-1$, then  $w_{\max}(n,k) = \left(\frac{k(n+1)}{2}\right)^n$.
If $n$ is even and $k$ is odd such that  $k \geq n-1$, then
$w_{\max}(n,k) = \left(\frac{k^2(n+1)^2 -1}{4}\right)^{n/2}$.
\end{theorem}

The corresponding Python function for computing $w_{\max}(n,k)$ (for $k\geq 2$):

\begin{verbatim}
from itertools import permutations, combinations_with_replacement
def wmax(n,k): # compute wmax(n,k) for k > 1
    ntuple, count = tuple(range(1,n+1)), 0
    for s in combinations_with_replacement(permutations(ntuple,n),k-2):
        t = list(ntuple)
        for d in s:
            for i in range(n):
                t[i] += d[i]
        t.sort()
        w = 1
        for i in range(n):
            w *= (n-i)+t[i]
        if w > count:
            count = w
    return count
\end{verbatim}

Partial list of values of $w_{\max}(n,k)$ are listed in Table \ref{tbl:wmax}. 
The antidiagonals of Table \ref{tbl:wmax} can be found in \url{https://oeis.org/A331988} (OEIS sequence A331988).

\begin{table}[htbp]
\begin{center}
\tiny
\begin{tabular}{r||r|r|r|r|r|r|r|r|r|r|r|r|}
& $k=1$ & 2 & 3&4&5&6&7&8&9&10&11&12\\
\hline\hline
$n=1$&1    & 2&     3  &   4 &    5  &   6 &    7  &   8   &  9 &    10  &   11   &  12   \\ \hline
2& 2  &   9 &    20 &    36   & 56  &  81 &   110  &  144 &   182  &  225 &   272 &  324 \\ \hline
3 &6 &   64  &  216 &   512  &  1000 &  1728 &  2744  & 4096 &  5832 & 8000 & 10648 & 13824   \\ \hline
4&24&    625   & 3136 &   10000  & 24336 &  50625 & 93636 & 160000 & 256036 &390625 &571536 &810000   \\ \hline
5 & 120   & 7776  &  59049 &  248832 &  759375 & 1889568 & 4084101  &7962624&14348907&24300000&39135393&60466176\\\hline
6 & 720  &  117649  & 1331000 &  7529536 & 28652616 &85766121 &216000000&481890304&976191488&1838265625&3254952168&5489031744\\\hline
\end{tabular}
\end{center}

\begin{center}
\tiny
\begin{tabular}{r||r|r|r|r|r|r|r|r|r|}
& $k=1$ & 2 & 3&4&5&6&7&8\\
\hline\hline
n = 7 & 5040 &  2097152 &  35831808  &268435456  &  1280000000 &4586471424 &13492928512&34359738368 \\\hline
8 & 40320 & 43046721 & 1097199376 &11019960576 & &282429536481&968381956096&2821109907456\\\hline
9 & 362880 &1000000000 & 38443359375 &512000000000 & &19683000000000&&262144000000000 \\\hline
10 & 3628800 & 25937424601 & 1488827973632 & 26559922791424& &1531578985264449&&27197360938418176\\\hline
11 & 39916800 & 743008370688 & 64268410079232 & 1521681143169024& &131621703842267136&&3116402981210161152\\\hline
\end{tabular}
\end{center}

\caption{Partial list of $w_{\max}(n,k)$.}\label{tbl:wmax}
\end{table}

\subsection{Permutation sets achieving $w_{\max}(n,k)$}
The next sections list, for each $n$ and $k$, a $k$-set of permutations (there may be many) that achieves $w_{\max}(n,k)$.  
We omit the cases where $n\leq2$, $k\leq 2$, even $k$ or $k \equiv 0 \mod n$ as they were discussed above.

\subsubsection{$k=3$}
\begin{itemize}
\item $n=4$: $(1234, 2341, 4213)$
\item $n=5$: $(12345, 34512, 53142)$
\item $n=6$: $(123456, 345612, 642153)$
\item $n=7$: $(1234567, 4567123, 7531642)$
\item $n=8$: $(12345678, 45378612, 86731254)$
\item $n=9$: $(123456789, 564897231, 978312645)$
\item $n=10$: $(123456789a, 56479a8123, a895312764)$
\item $n=11$: $(123456789ab, 6758ab92413, b9a63128574)$
\end{itemize}

\subsubsection{$k=5$}
\begin{itemize}
\item $n=3$: $(123, 123, 231, 312, 321)$
\item $n=4$: $(1234, 1234, 2341, 4213, 4321)$
\item $n=6$: $(123456, 123456, 342561, 645213, 654312)$
\item $n=7$: $(1234567, 1234567, 4536712, 7562143, 7654231)$
\end{itemize}

\subsubsection{$k=7$}
\begin{itemize}
\item $n=3$: $(123, 123, 123, 231, 312, 321, 321)$
\item $n=4$: $(1234, 1234, 1234, 2341, 4213, 4321, 4321)$
\item $n=5$: $(12345, 12345, 12345, 34251, 53412, 54312, 54321)$
\item $n=6$: $(123456, 123456, 123456, 342561, 645213, 654312, 654321)$
\end{itemize}

\subsubsection{$k=9$}
\begin{itemize}
\item $n=4$: $(1234, 1234, 1234, 1234, 2341, 4213, 4321, 4321, 4321)$
\item $n=5$: $(12345, 12345, 12345, 12345, 34251, 53412, 54312, 54321, 54321)$
\end{itemize}

\subsubsection{$k=11$}
\begin{itemize}
\item $n=3$: $(123, 123, 123, 123, 123, 231, 312, 321, 321, 321, 321)$
\item $n=4$: $(1234, 1234, 1234, 1234, 1234, 2341, 4213, 4321, 4321, 4321, 4321)$
\end{itemize}

\subsubsection{$k=13$}
\begin{itemize}
\item $n=3$: $(123, 123, 123, 123, 123, 123, 231, 312, 321, 321, 321, 321, 321)$
\item $n=4$: $(1234, 1234, 1234, 1234, 1234, 1234, 2341, 4213, 4321, 4321, 4321, 4321, 4321)$
\end{itemize}

\subsubsection{$k=15$}
\begin{itemize}
\item $n=4$: ($1234$, $1234$, $1234$, $1234$, $1234$, $1234$, $1234$, $2341$, $4213$, $4321$, $4321$, $4321$, $4321$, $4321$, $4321$)
\end{itemize}

\subsubsection{$k=17$}
\begin{itemize}
\item $n=3$: ($123$, $123$, $123$, $123$, $123$, $123$, $123$, $123$, $231$, $312$, $321$, $321$, $321$, $321$, $321$, $321$, $321$)
\item $n=4$: ($1234$, $1234$, $1234$, $1234$, $1234$, $1234$, $1234$, $1234$, $2341$, $4213$, $4321$, $4321$, $4321$, $4321$, $4321$, $4321$, $4321$)
\end{itemize}

\end{document}